\documentclass[11pt,a4paper]{article}

\usepackage{latexsym}
\usepackage{amssymb,amsfonts,amscd,amsmath}
\usepackage{epsfig}

\newcommand{\bC}{{\mathbb C}}
\newcommand\sD{{\mathcal D}}
\newcommand\sV{{\mathcal V}}

\begin{document}

\title{Solution of polynomial systems derived from differential equations}

\author{\textbf{E. L. Allgower}, Fort Collins, 
\textbf{D. J. Bates}, Notre Dame,\\ 
\textbf{A. J. Sommese}, Notre Dame, and 
\textbf{C. W. Wampler}, Warren}

\maketitle

\begin{abstract}
Nonlinear two-point boundary value problems arise in numerous areas
of application. The existence and number of solutions for various
cases has been studied from a theoretical standpoint.  These results
generally rely upon growth conditions of the nonlinearity. However,
in general, one cannot forecast how many solutions a boundary value
problem may possess or even determine the existence of a solution.
In recent years numerical continuation methods have been developed
which permit the numerical approximation of all complex solutions of
systems of polynomial equations. In this paper, numerical
continuation methods are adapted to numerically calculate the
solutions of finite difference discretizations of nonlinear
two-point boundary value problems.  The approach taken here is to
perform a homotopy deformation to successively refine
discretizations.  In this way additional new solutions on finer
meshes are obtained from solutions on coarser meshes.  The
complicating issue which the complex polynomial system setting
introduces is that the number of solutions grows with the number of
mesh points of the discretization.  To counter this, the use of
filters to limit the number of paths to be followed at each stage is
considered.
\end{abstract}



\section{Introduction}
Consider a two-point boundary value problem on the interval
$[a,b]\subset \mathbb{R}$,
\begin{equation}\label{bvp1}
y''=f(x,y,y'),
\end{equation}
with boundary conditions $y(a)=\alpha$ and $y(b)=\beta$.  The
standard central difference approximation with a uniform mesh may be
used to approximate solutions to (\ref{bvp1}).  In particular, let
$N$ be a positive integer, $h:=\frac{b-a}{N+1}$, and $x_i:=a+ih$ for
$i=0, \dots, N+1$. Setting $y_0=\alpha$ and $y_{N+1}=\beta$, the
discretization of (\ref{bvp1}) takes the form of the system $\sD_N$:
\[
\begin{array}{llll}
y_0&-\hskip.15in 2y_1&+\hskip.15in y_2 &= h^2f(x_1,y_1,\frac{y_2-y_0}{2h})\\
\hskip.05in \vdots &\hskip.35in \vdots &\hskip.3in \vdots & = \hskip.45in \vdots\\
y_{N-1}&-\hskip.15in 2y_N&+\hskip.15in y_{N+1} & = h^2f(x_N,y_N,\frac{y_{N+1}-y_{N-1}}{2h})
\end{array}
\]
A solution $y(x)$ of (\ref{bvp1}) may then be approximated by an $N$-tuple of real
numbers $(y_1,\dots, y_N)$ such that $y_i\approx y(x_i)$ for $i=1, \dots, N$.

Depending upon the nonlinearity $f$, equation (\ref{bvp1}) may have
no solutions, a unique solution, multiple solutions, or even
infinitely many solutions. There are many existence theorems for
solutions of such equations subject to growth conditions on $f$, but
even when existence is known, the number of solutions often is not.
Furthermore, a discretization such as $\sD_N$ may have spurious
solutions that do not converge to a solution to (\ref{bvp1}) as
$N\to\infty$. On the other hand, if $f$ is sufficiently smooth, a
solution $y$ to (\ref{bvp1}) is eventually approximated  with
$\mathcal{O}(h^2)$ accuracy on the mesh by some solution $\bar{y}\in
\mathbb{R}^N$.

The purpose of the present paper is to give a relatively secure
numerical technique for finding the solutions of a general class of
two-point boundary value problems without requiring highly refined
meshes. The technique involves performing successive homotopy
deformations between discretizations with increasingly many mesh
points, as suggested in \cite{A2}. By restricting our attention to
problems having polynomial nonlinearity, including the case of a
polynomial approximation to a smooth nonlinearity, we can often
assure that all solutions are found at each stage of the algorithm.
Even when we do not guarantee all solutions, our method generates
multiple solutions that in test cases include approximations to all
known solutions. While Gr\"obner basis methods (see \cite{CLO}) or
cellular exclusion methods (see \cite{G}) could be applied to solve
the polynomial discretizations, we chose to use homotopy
continuation due to its ability to handle polynomial systems in many
variables and the ease with which it allows us to generate solutions
on a refined mesh from the solutions on the previous mesh.  Although 
other numerical methods treating two-point boundary value problems have 
been developed (see \cite{K} and \cite{RS}),
such methods require satisfactory initial solution estimates.  The present
technique provides such initial estimates.

Here is a sketch of our bootstrapping process, which will be discussed 
in more detail in the subsequent section:

\begin{enumerate}
\item Find all solutions of the discretization $\sD_N$ for some
small $N$. The size of $N$ needs only to be large enough that the
discretization is consistent; it could be as small as $N=1$.

\item Discard all unreasonable solutions, e.g., solutions
which do not possess properties which exact solutions may be known
to have. Let us denote the set of solutions which are kept by
$\sV_N$.

\item If the
mesh size is not yet sufficiently small or the cardinality of
$\sV_N$ has not yet stabilized, add a mesh point to obtain the
discretization $\sD_{N+1}$.  Use the solutions in $\sV_N$ to
generate solutions $\sV_{N+1}$ of $\sD_{N+1}$ and then return to
Step 2.

\item Once the mesh size $h$ is sufficiently small and the cardinality of
$\sV_N$ becomes stable, refine the solutions with a fast nonlinear
solver, using starting values obtained by interpolating the
solutions in $\sV_N$.
\end{enumerate}
This paper focuses primarily on the implementation of Step 3 of the
above scheme.  In particular, we consider the homotopy function
\begin{align*}
  H_{N+1}&(y_1,\ldots,y_{N+1},t):=\\
   &\left[
   \begin{array}{rcl}
   y_0-2y_1+y_2&-&h(t)^2f\left(x_1(t),y_1,\frac{y_2-y_0}{2h(t)}\right)\\
   &\vdots&\\
   y_{N-2}-2y_{N-1}+y_{N}&-&h(t)^2f\left(x_{N-1}(t),y_{N-1},\frac{y_{N}-y_{N-2}}{2h(t)}\right)\\
   y_{N-1}-2y_{N}+Y_{N+1}(t)&-&h(t)^2f\left(x_{N}(t),y_{N},\frac{Y_{N+1}(t)-y_{N-1}}{2h(t)}\right)\\
   y_{N}-2y_{N+1}+Y_{N+2}(t)&-&h(t)^2f\left(x_{N+1}(t),y_{N+1},\frac{Y_{N+2}(t)-y_{N}}{2h(t)}\right)\\
   \end{array}
   \right]
\end{align*}
with
 $$
  \begin{array}{rcl}
   y_0&:=&\alpha\\
   h(t)&:=&t\left(\frac{b-a}{N+1}\right)+(1-t)\left(\frac{b-a}{N+2}\right)\\
   Y_{N+1}(t)&:=&(1-t)y_{N+1}+\beta t\\
   Y_{N+2}(t)&:=&\beta (1-t)\\
   x_i(t)&:=&a+ih(t),\qquad i=1,\ldots,N+1.
   \end{array}
 $$
At $t=0$ this is the system $\sD_{N+1}$.  At $t=1$, it can be
interpreted as the system $\sD_N$ with a new mesh point having the
value $y_{N+1}$ at $x=b$ and a new right-hand boundary at $x=b+h(1)$
having value $Y_{N+2}(1)$. The incompatibility of the old boundary
condition at $x=b$ and the new one at $x=b+h(1)$ is accommodated by
the presence of both $y_{N+1}$ and $Y_{N+1}$, which are not
necessarily equal. As $t$ goes from 1 to 0, the mesh points are
squeezed back inside the interval $[a,b]$, and the right-hand
boundary condition $y(b)=\beta$ is transferred from $Y_{N+1}$ to
$Y_{N+2}$ as $Y_{N+1}$ is forced to equal $y_{N+1}$.

To find solutions of $\sD_{N+1}$, we use continuation to track
solutions of $H_{N+1}$ as $t$ goes from 1 to 0.  At $t=1$, we have a 
list $\sV_N$ of solutions $(y_1,\ldots,y_N)$ satisfying the first
$N$ equations of $H_{N+1}$, while the final equation is
 $$
  y_{N}-2y_{N+1}=h(t)^2f\left(b,y_{N+1},\frac{-y_{N}}{2h}\right),
 $$
which is the only place where $y_{N+1}$ appears.  For each solution
of $(y_1,\ldots,y_N)$ in $\sV_N$, we may use this equation to find
corresponding solution values for $y_{N+1}$.  These are the start
points of continuation paths leading to solutions of $\sD_{N+1}$.

The framework above does not change in any of its essentials if we
prescribe a different function for $Y_{N+2}(t)$.  For example, the
constant function $Y_{N+2}(t):=\beta$ was used for all examples 
below.  Although other alternatives are 
theoretically feasible, none were tested.  The essential feature of 
$Y_{N+2}(t)$ is that it goes to $\beta$ as $t$ goes to $0$.

By the implicit function theorem, a nonsingular solution $y=y^*$ to
$H_{N+1}(y,1)=0$ will continue uniquely in the neighborhood of $t=1$
to a nonsingular solution path $y(t)$ satisfying $H_{N+1}(y(t),t)=0$
with $y(1)=y^*$. This does not mean, however, that the path remains
nonsingular all the way to $t=0$, which is what we require to follow
the path reliably with numerical continuation.  To skirt this
difficulty, as discussed in Chapter 7 of \cite{SoWa}, it is
sufficient to insert a random $\gamma\in \bC$ into the homotopy to
obtain the variant
\begin{align}
  &H_{N+1}(y_1,\ldots,y_{N+1},t):=\nonumber\\
   &\left[
   \begin{array}{rcl}
   \Gamma(t)\left(y_0-2y_1+y_2\right)&-&h(t)^2f\left(x_1(t),y_1,\frac{y_2-y_0}{2h(t)}\right)\\
   &\vdots&\\
   \Gamma(t)\left(y_{N-2}-2y_{N-1}+y_{N}\right)&-&h(t)^2f\left(x_{N-1}(t),y_{N-1},\frac{y_{N}-y_{N-2}}{2h(t)}\right)\\
   \Gamma(t)\left(y_{N-1}-2y_{N}\right)+Y_{N+1}(t)&-&h(t)^2f\left(x_{N}(t),y_{N},\frac{Y_{N+1}(t)-y_{N-1}}{2h(t)}\right)\\
   \Gamma(t)\left(y_{N}-2y_{N+1}+\beta\right)&-&h(t)^2f\left(x_{N+1}(t),y_{N+1},\frac{\beta-y_{N}}{2h(t)}\right)\\
   \end{array}
   \right] \label{Eq:HN1}
\end{align}
with
 $$
  \begin{array}{rcl}
   \Gamma(t)&:=&\gamma ^2 t+(1-t)\\
   h(t)&:=&\gamma
   t\left(\frac{b-a}{N+1}\right)+(1-t)\left(\frac{b-a}{N+2}\right)\\
   Y_{N+1}(t)&:=&(1-t)y_{N+1}+\gamma^2\beta t\\
   x_i(t)&:=&a+ih(t),\qquad i=1,\ldots,N+1.
   \end{array}
 $$

The work of the second author was supported by the National Science Foundation under Grant No.\ 0105653 and
Grant No.\ 0410047; and a fellowship from the Arthur J. Schmitt Foundation.  The work of the third author
was supported by the National Science Foundation under Grant No.\ 0105653 and Grant No.\ 0410047; and
the Duncan Chair of the University of Notre Dame.  The work of the fourth author was supported by 
the National Science Foundation under Grant No.\ 0410047.

\section{The case of polynomial nonlinearity} Let's specialize the homotopy
of equation (\ref{Eq:HN1}) to the case when $f(x,y,y')$ is a real
polynomial $p(y)$.  Then, the right-most term of the $i^{\rm th}$
entry in $H_{N+1}(y,t)$ becomes just $h^2(t)p(y_i)$. This
restriction to the polynomial case allows us to conveniently obtain
the start points for $H_{N+1}(y_1,\ldots,y_{N+1},1)=0$ by solving
the polynomial
$$y_{N}-2y_{N+1}+\beta-\left(\frac{b-a}{N+1}\right)^2p(y_{N+1})=0$$
for $y_{N+1}$ given $y_N$ from the solutions in $\sV_N$.

Let $d=\deg p(y)$.  We see that, in general, over the complex
numbers, we will obtain $d$ values of $y_{N+1}$ for every point in
$\sV_N$.  Suppose that at each stage of the algorithm these all
continue to finite, nonsingular solutions of $\sD_{N+1}$.  Then, the
solution list $\sV_N$ will have $d^N$ entries.  While this gives an
exhaustive enumeration of the solutions of the discretized problem,
the exponential growth in the length of the solution list cannot be
practically sustained as $N$ increases. However, it is often the
case that most of the solutions at a given stage do not exhibit
various properties required of solutions to the two-point boundary
value problem at hand, leading to filtering rules.  Depending upon
the problem at hand, there are a variety of filtering rules that may
be implemented to determine which solutions in $\sV_N$ may be
discarded as start solutions for the subsequent homotopy.

For small $N$ we can contemplate retaining all solutions.  It is
reasonable to ask whether the above procedure is guaranteed to
generate all solutions of the discretized system.  The answer, in
general, is no, but we can say that if $\sV_N$ has $d^N$ distinct,
nonsingular solutions, then it is clear that all solutions have been
found, as B\'ezout's theorem states that this is the greatest number
possible.  Indeed, for our test problems, we have found that this
behavior is typical.

For larger $N$, a filter becomes necessary. One that is always
available is to take the discretization of
 the derivative of $y'''=p'(y)y'$ of $y''=p(y)$, and throw away
 $y\in\sV_N$ for which this is large.  To get the discretization
 we could use the central difference approximations
  $$
  y'(x_i)=\frac{y_{i+1}-y_{i-1}}{2h},
  $$
  and
  $$
  y'''(x_i)=\frac{y_{i+2}-2y_{i+1}+2y_{i-1}-y_{i-2}}{2h^3}
  $$
applied only at the mesh points $y_2,\ldots,y_{N-1}$.  So we would
throw away the point $z=(y_1,\ldots,y_N)\in\sV_N$ if
  $$
  \sum_{i=2}^{N-1}\left|\frac{y_{i+2}-2y_{i+1}+2y_{i-1}-y_{i-2}}{2h^3}
     -p'(y_i)\frac{y_{i+1}-y_{i-1}}{2h}\right| >\epsilon_2
  $$
for some $\epsilon_2>0$.  Naturally, one drawback to such a filter is the need to specify $\epsilon_2$.

Other filters may be derived from known properties of the solutions
of the problem at hand.  For example, it may be known that solutions
are symmetric about $x = \frac{a+b}{2}$, are always positive,
oscillate with a specific period, or exhibit some other
easily-detected behavior.  For example, a filter based on symmetry
is considered in Section 3.2.  Although one may be tempted to
discard solutions having nonzero complex part, this is not a valid
filtering rule. The problem in Section 3.3 below has non-real
solutions in $\sV_N$ that are tracked to real solutions in
$\sV_{N+1}$. Similarly, it is possible that oscillating solutions
may arise from a sequence of non-oscillating solutions and that
similar problems may occur with other filters.  Thus, the use of
filters may be computationally beneficial, but with it comes the
risk of not finding all real solutions to the problem.

Thus we have the final version of the algorithm:

\begin{quote}
{\bf Algorithm 1}
\begin{enumerate}
\item For $N = 1$, $H_1(y_1, 1)$ is a single polynomial in $y_1$,
which may be solved with any one-variable method to produce $\sV_1$.

\item For $N = 2, 3, ...$, until some desired behavior has occurred:

\begin{enumerate}
\item Form the homotopy $H_{N}(y_1,\ldots,y_{N},t)$.

\item Solve the last polynomial of $H_{N}(y_1,\ldots,y_{N},t)$ for
$y_{N}$ using each solution in $\sV_{N-1}$, thereby forming the set
$S$ of the start solutions for $H_{N}(y_1,\ldots,y_{N},t)$.

\item Track all paths beginning at points in $S$ at $t=1$.  The
set of endpoints of these paths is $\sV_N$.

\item If desired, apply a filter to $\sV_N$ to reduce the number of
paths to be tracked in stage $N+1$.
\end{enumerate}
\item Refine the solutions with a nonlinear solver, if desired.
\end{enumerate}

\end{quote}

\section{Numerical experiments}

The following experiments were run using Bertini, a software package
under development by the last three authors for the study of
numerical algebraic geometry.  Although Bertini was written to make
use of multiprecision adaptively, each of the following experiments
ran successfully using only 16 digits of precision.

In the following, $N$ denotes the number of mesh points, SOLS($N$)
denotes the total number of solutions (real or complex), and
REAL($N$) denotes the number of real solutions.  For $N > 1$, the
number of paths tracked from stage $N-1$ is $d\cdot{\rm SOLS}(N-1)$.
A solution is considered to be real if the imaginary part at each
mesh point is zero to at least eight digits.

\subsection{A basic example}

As a first example, consider the following two-point boundary value
problem
\begin{equation}\label{bvp2}
y'' = 2y^3
\end{equation}
with boundary conditions $y(0) = \frac{1}{2}$ and $y(1) = \frac{1}{3}$.

There is a unique solution, $y = \frac{1}{x+2}$, to (\ref{bvp2}).
Our method produces one real solution among a total of $3^N$
solutions found for $N = 1, \dots, 9$.  Furthermore, the error
between the computed solution and the unique exact solution behaves
as $\mathcal{O}(h^2)$.  Refer to Table 1 for details.

\begin{center}
\begin{tabular}{|c|c|c|c|}
\hline
$N$&Maximal error at any mesh point&$h^2$&Maximal error/${h^2}$\\
\hline
3& 1.570846e-04& 4.000000e-02& 3.927115e-03\\
4& 1.042635e-04& 2.777778e-02& 3.753486e-03\\
5& 7.069710e-05& 2.040816e-02& 3.464158e-03\\
6& 5.348790e-05& 1.562500e-02& 3.423226e-03\\
7& 4.078910e-05& 1.234568e-02& 3.303917e-03\\
8& 3.230130e-05& 1.000000e-02& 3.230130e-03\\
9& 2.624560e-05& 8.264463e-03& 3.175718e-03\\
\hline
\multicolumn{4}{c}{\small Table 1:  Evidence of $\mathcal{O}(h^2)$ convergence for Problem (\ref{bvp2}).}\\
\end{tabular}
\end{center}

\subsection{A more sophisticated example}

Consider the problem
\begin{equation}\label{bvp3}
y'' = -\lambda\left(1+y^2\right)
\end{equation}
with zero boundary conditions, $y(0) = 0$ and $y(1) = 0$, and $\lambda > 0$.

According to \cite{La}, any solutions to this problem must be
symmetric about $x=\frac{1}{2}$, so we have a special filter.
Furthermore, it is known that there are two solutions if $\lambda <
4$, a unique solution if $\lambda = 4$, and no solutions if $\lambda
> 4$.  Without using a filter, the expected number of real solutions
in the first and last cases were confirmed computationally (for
$\lambda = 2$ and $\lambda = 6$ with $N \leq 17$), and the computed
solutions were symmetric as anticipated. From B\'ezout's theorem,
one would expect to find at most $2^N$ complex solutions at each
stage $N$, and this is precisely the total number of complex
solutions found.  When $\lambda \approx 4$, the Jacobian of the
associated polynomial system is rank-deficient, so regular
path-tracking techniques fail.

Tracking all $2^{17}$ paths for $N=17$ took just under an hour of CPU time on a single processor Pentium 4, 3 GHz machine running Linux.
At this rate, ignoring the time-consuming data management part of the algorithm, it would take well over one year to track all
$2^{30}$ paths for $N=30$ mesh points.  As discussed in Section 2, filtering rules may be used to dramatically reduce the number
of paths to be tracked at each stage.  A filter forcing $\left|\left|y_1\right| - \left|y_N\right|\right| < 10^{-8}$ was applied to the
case $\lambda = 2$.  This cut the path-tracking time to less than half a second for $N=17$ mesh points.  This drastic reduction
in time for path-tracking as well as data management allowed for the confirmation of the existence of two real solutions for up to $100$ mesh points.
Despite the size of the polynomial system when $N = 100$, each path took less than $4$ seconds to track from $t=1$ to $t=0$.  A graph of the
two real solutions for $N=20$ mesh points is given in Figure 1.

\begin{figure}
\centerline{
\epsfxsize=3in
\epsfysize=3in
\epsfbox{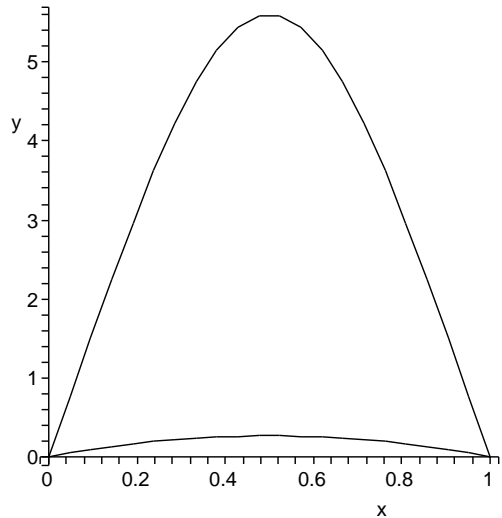}
}
\caption{The real solutions of (\ref{bvp3}) with $N=20$.}
\end{figure}

\subsection{A problem with infinitely many solutions}

It was shown in \cite{Co} that the two-point boundary value problem
\begin{equation}\label{bvp4}
y'' = -\lambda y^3, \qquad y(0)=y(1)=0,
\end{equation}

\noindent with $\lambda > 0$ has infinitely many oscillating real
solutions on the interval $[0,1]$.  Moreover, the solutions occur in
pairs in the sense that $-y$ is a solution whenever $y$ is a
solution. Hence, together with the trivial solution $y=0$, we expect
always to have an odd number of solutions.  That was confirmed
computationally, as shown in Table 2. Only the case of $\lambda=1$
is displayed as all other cases are identical modulo scaling. It may
be observed that the number of real solutions found by Bertini grows
without bound for this problem, as the number of mesh points
increases.  In fact, beyond some small value of $N$, the number of
real solutions approximately doubles for each subsequent value of
$N$.

\begin{center}
\begin{tabular}{|c|c|c|}
\hline
$N$&SOLS($N$)&REAL($N$)\\
\hline
1&3&3\\
2&3&3\\
3&9&3\\
4&27&7\\
5&81&11\\
6&243&23\\
7&729&47\\
8&2187&91\\
\hline
\multicolumn{3}{c}{\small Table 2:  Solutions of (\ref{bvp4})}\\
\end{tabular}
\end{center}

\subsection{The Duffing problem}

One representation (see \cite{D}) of the Duffing problem is the two-point 
boundary value problem
\begin{equation}\label{bvp5}
y'' = -\lambda \sin\left(y\right)
\end{equation}
on the interval $\left[0,1\right]$ with  $y(0) = 0$, $y(1) = 0$, and
$\lambda > 0$.  Since our attention is restricted to polynomial
nonlinearity only, we approximate $\sin(y)$ by truncating its power
series expansion, yielding the problem
\begin{equation}\label{bvp6}
y'' = -\lambda \left(y - \frac{y^3}{6}\right)
\end{equation}
using two terms or
\begin{equation}\label{bvp7}
y'' = -\lambda \left(y - \frac{y^3}{6} + \frac{y^5}{120}\right)
\end{equation}
using three terms.

It is known that there are $2k+1$ real solutions to the exact
Duffing problem  (\ref{bvp5}) when $k\pi < \lambda <
\left(k+1\right)\pi$. For a given value of $\lambda$, the $2k+1$
real solutions include the trivial solution $y \equiv 0$ and $k$
pairs of solutions $\left(y_1(x), y_2(x)\right)$ such that $y_1(x) =
-y_2(x)$.  Each pair oscillates with a different period.  As two-
and three-term Taylor series truncations for $\sin(y)$ do not
approximate $\sin(y)$ well outside of a small neighborhood, the
solutions to (\ref{bvp6}) and (\ref{bvp7}) may behave quite
differently than those of (\ref{bvp5}).

Table 3 indicates the number of real solutions found for problems
(\ref{bvp6})  and (\ref{bvp7}) for $\lambda = 0.5\pi$, $1.5\pi$, and
$2.5\pi$.  All solutions have either odd or even symmetry about $x=\frac{1}{2}$,
so we again used the filter $\left|\left|y_1\right| - \left|y_N\right|\right| < 10^{-8}$.  
The filter was first applied when $N = 4$, so the number of real solutions reported in
each case of Table 3 is the number of real solutions found for $N
\geq 5$.  For $\lambda = 0.5\pi$ and $\lambda = 1.5\pi$, there were
more real solutions found for (\ref{bvp7}) than predicted for the
exact problem (\ref{bvp5}).  However, the computed solutions in each
case included one pair of solutions that oscillated wildly.  These
poorly-behaved solutions are readily identified by the $y'''$ filter
discussed in Section 2: for $N = 25$ mesh points, they had residuals
four orders of magnitude larger than those of the well-behaved
solutions.

\begin{center}
\begin{tabular}{|c|c|c|c|}
\hline
$\lambda$& $f(y) = y-\frac{y^3}{6}$ & $f(y) = y - \frac{y^3}{6} + \frac{y^5}{120} $ & $ f(y) = \sin(y)$\\
\hline
$0.5\pi$& 1 & 3 & 1 \\
$1.5\pi$& 1 & 5 & 3 \\
$2.5\pi$& 1 & 5 & 5 \\
\hline
\multicolumn{4}{c}{\small Table 3:  Number of real solutions for approximations of the Duffing problem.}\\
\end{tabular}
\end{center}

\subsection{The Bratu problem}
The Bratu problem on the interval
$\left[0,1\right]$ has the form
\begin{equation}\label{bvp8}
y'' = -\lambda e^y,\qquad y(0)=y(1)=0,
\end{equation}
with $\lambda > 0$.  As in the case of the Duffing problem, we make
the right-hand side polynomial by truncating the power series
expansion of $e^y$, yielding
\begin{equation}\label{bvp9}
y'' = -\lambda \left(1 + y + \frac{y^2}{2}\right)
\end{equation}

As discussed in \cite{D}, there are two real solutions if $\lambda$
is near zero and no real solutions if $\lambda$ is large. The real
solutions for small $\lambda$ are symmetric and nonnegative.  The
expected number and properties of the real solutions in the cases of
$\lambda = 0.5$ and $\lambda = 10$ were confirmed, and, as anticipated, 
$2^N$ total solutions were found in each case for $N=1, \dots, 15$.

\section{Discussion}

A new algorithm for finding the real solutions of a two-point
boundary value problem has been presented, and several examples have
been documented under the assumption of polynomial nonlinearity.
Furthermore, the use of filtering rules to drastically reduce the
computational work has been considered.  In each example presented,
the number of real solutions predicted by theory has been confirmed
computationally, although it was seen that the use of filters may
effect the number of real solutions discovered.

There are several variations to the algorithm that could be
considered in the future.  A more detailed analysis of the benefits
and drawbacks of the use of filters could be made.  Also, it is
possible to add extra mesh points at the left-hand end or middle of
the interval rather than the right.  Similarly, one new mesh point
could be added to each end simultaneously, yielding $\left(\deg
p(y)\right)^2$ starting solutions for each solution from the
previous stage. For that matter, non-uniform grids could be analyzed
with only mild changes to the formulation.  A similar algorithm
could also be developed for systems of differential equations.

\bigskip

\noindent Eugene L. Allgower

\noindent Department of Mathematics

\noindent Colorado Sate University

\noindent Fort Collins, CO 80523-1874

\noindent USA

\noindent allgower@math.colostate.edu

\bigskip

\noindent Daniel J. Bates

\noindent Department of Mathematics

\noindent University of Notre Dame

\noindent Notre Dame, IN 46556-4618

\noindent USA

\noindent dbates1@nd.edu

\bigskip

\noindent Andrew J. Sommese

\noindent Department of Mathematics
                                                                                                                                                                             
\noindent University of Notre Dame
                                                                                                                                                                             
\noindent Notre Dame, IN 46556-4618
                                                                                                                                                                             
\noindent USA

\noindent sommese@nd.edu

\bigskip

\noindent Charles W. Wampler

\noindent General Motors Research and Development

\noindent Mail Code 480-106-359

\noindent 30500 Mound Road

\noindent Warren, MI 48090-9055

\noindent USA

\noindent Charles.W.Wampler@gm.com

\end{document}